\documentclass[12pt]{article}
\usepackage{amsmath,amsthm,amsfonts,amssymb,amscd}
\def\qed{{\unskip\nobreak\hfil\penalty50
\hskip2em\hbox{}\nobreak\hfil$\square$
\parfillskip=0pt \finalhyphendemerits=0\par}\medskip}
\def\proof{\trivlist \item[\hskip \labelsep{\bf Proof\ }]}
\def\endproof{\null\hfill\qed\endtrivlist}
\def\lan{\langle}
\def\ran{\rangle}
\def\Ad{{\mathrm {Ad}}}
\def\Vir{{\mathrm {Vir}}}

\def\End{{\mathrm {End}}}
\def\Hom{{\mathrm {Hom}}}
\def\id{{\mathrm {id}}}

\def\Vir{{\mathrm {Vir}}}

\def\a{\alpha}

\def\e{\varepsilon}

\def\la{\lambda}
\def\La{\Lambda}
\def\phi{\varphi}

\def\Om{\Omega}

\def\r{{\rho}}

\def\emptyset{\varnothing}
\def\setminus{\smallsetminus}

\def\Diff{{\mathrm {Diff}}}
\def\Mob{{\rm\textsf{M\"ob}}}

\def\lan{\langle}
\def\ran{\rangle}
\def\Ad{{\mathrm {Ad}}}

\def\End{{\mathrm {End}}}
\def\Hom{{\mathrm {Hom}}}
\def\id{{\mathrm {id}}}

\def\Vir{{\mathrm {Vir}}}

\def\a{\alpha}

\def\e{\varepsilon}

\def\l{\lambda}
\def\La{\Lambda}
\def\phi{\varphi}

\def\Om{\Omega}

\def\r{{\rho}}

\newtheorem{theorem}{Theorem}[section]
\newtheorem{lemma}[theorem]{Lemma}
\newtheorem{conjecture}[theorem]{Conjecture}
\newtheorem{corollary}[theorem]{Corollary}
\newtheorem{definition}[theorem]{Definition}

\newtheorem{proposition}[theorem]{Proposition}
\newtheorem{remark}[theorem]{Remark}

\def\col{{\mathrm{col}}}
\def\emptyset{\varnothing}
\def\setminus{\smallsetminus}

\def\exp{{\mathrm {exp}}}

\def\Diff{{\mathrm {Diff}}}
\def\Mob{{\rm\textsf{M\"ob}}}

\def\res{\!\restriction\!}
\def\A{{\cal A}}

\def\B{{\cal B}}
\def\C{{\cal C}}
\def\D{{\cal D}}

\def\I{{\cal I}}

\def\L{{\cal L}}

\def\H{{\cal H}}

\def\Z{{\mathbb Z}}

\renewcommand{\qed}{\ \hfill $\blacksquare$}

\newcommand{\bdefin}{\begin{definition}}
\newcommand{\blemma}{\begin{lemma}}
\newcommand{\bprop}{\begin{proposition}}
\newcommand{\btheor}{\begin{theorem}}
\newcommand{\bcoro}{\begin{corollary}}
\newcommand{\bconj}{\begin{conjecture}}
\newcommand{\edefin}{\end{definition}}
\newcommand{\elemma}{\end{lemma}}
\newcommand{\eprop}{\end{proposition}}
\newcommand{\etheor}{\end{theorem}}
\newcommand{\ecoro}{\end{corollary}}
\newcommand{\econj}{\end{conjecture}}
\newcommand{\brem}{\begin{remark}}
\newcommand{\erem}{\end{remark}}

\newcommand{\ba}{\begin{array}}
\newcommand{\ea}{\end{array}}
\newcommand{\bea}{\begin{eqnarray}}
\newcommand{\eea}{\end{eqnarray}}
\newcommand{\bean}{\begin{eqnarray*}}
\newcommand{\eean}{\end{eqnarray*}}

\renewcommand{\mod}{\mbox{mod}}

\title{\huge An application of Mirror extensions \\}
\author{
{\sc Feng Xu}\footnote{Supported in part by NSF.}\\
Department of Mathematics\\
University of California at Riverside\\
Riverside, CA 92521\\
E-mail: {\tt xufeng@math.ucr.edu}}
\begin{document}
\date{}
\maketitle

\begin{abstract}
In this paper we apply our previous results of mirror extensions to
obtain   realizations of three modular invariants constructed by A.
N. Schellekens by holomorphic conformal nets with central charge
equal to $24$.\par 2000MSC:81R15, 17B69.
\end{abstract}
%\maketitle
%\thanks
%{I'd like to thank Professor Vaughan Jones for very valuable discussions.
%This work is partially supported by NSF.}
%\endthanks
\newpage

\section{Introduction}
Partition functions of chiral rational conformal field theories
(RCFT) are modular invariant (cf. \cite{Z}). However there are
examples of ``spurious"  modular invariants which do not correspond
to any RCFT (cf.\cite{BE4}, \cite{SY} and \cite{FSS}) . It is
therefore an interesting question to decide which modular invariants
can be realized in RCFT. For many interesting modular invariants
this question was raised  for an example in \cite{Sch} and more
recently in \cite{EW}. For results on related questions, see
\cite{BE3}, \cite{BE4}, \cite{Reh1}, \cite{KL},\cite{KL2} and
\cite{KLPR} for a partial list.
\par In this paper we examine the holomorphic modular invariants
with central charge 24 constructed by A. N. Schellekens in
\cite{Sch}. Besides modular invariance, A. N. Schellekens showed
that his modular invariants passed an impressive list of checks from
tracial identities which strongly suggested that his modular
invariants can be realized in chiral RCFT. Some of Schellekens's
modular invariants were constructed using level-rank duality. In
\cite{Xm} we proved a general theorem on mirror extensions (cf. Th.
\ref{mainmirror}) which included modular invariants from level-rank
duality (cf. \S\ref{lr2}). It is therefore an interesting question
to see if mirror extensions can provide chiral RCFT realization of
some of Schellekens's modular invariants. Our main result in this
paper is to show that three of Schellekens's modular invariants can
be realized by holomorphic conformal nets (cf. Th. \ref{main}):
these nets are constructed by simple current extensions (cf.
\S\ref{simpleextension}) of mirror extensions. Our results strongly
suggest that there should be Vertex Operator Algebras which realize
these modular invariants. We expect our methods to apply to other
modular invariants in the literature, especially when level-rank
duality plays a role.\par This paper is organized as follows: after
a preliminary section on nets, mirror extensions and simple current
extensions,  we examine three of Schellekens's modular invariants in
\cite{Sch}, and obtain realization of these invariants as simple
current extensions of three mirror extensions. We end with two
conjectures about holomorphic conformal nets with central charge 24
which are motivated by \cite{FLM} and \cite{Sch}, and we hope that
these conjectures will stimulate further research.
\section{Preliminaries}
%\subsection{Conformal nets}

%\section{Preliminaries}
\subsection{Preliminaries on sectors}

Given an infinite factor $M$, the {\it sectors of $M$}  are given by
$$\text{Sect}(M) = \text{End}(M)/\text{Inn}(M),$$
namely $\text{Sect}(M)$ is the quotient of the semigroup of the
endomorphisms of $M$ modulo the equivalence relation: $\rho,\rho'\in
\text{End}(M),\, \rho\thicksim\rho'$ iff there is a unitary $u\in M$
such that $\rho'(x)=u\rho(x)u^*$ for all $x\in M$.

$\text{Sect}(M)$ is a $^*$-semiring (there are an addition, a product and
an involution $\rho\rightarrow \bar\rho$)
equivalent to the Connes correspondences (bimodules) on
$M$ up to unitary equivalence. If $\r$ is
an element of $\text{End}(M)$ we shall denote by $[\r]$
its class in $\text{Sect}(M)$. We define
$\text{Hom}(\r,\r')$ between the objects $\r,\r'\in \End(M)$
by
\[
\text{Hom}(\r,\r')\equiv\{a\in M: a\r(x)=\r'(x)a \ \forall x\in M\}.
\]
We use $\langle  \lambda , \mu \rangle$ to denote the dimension of
$\text{\rm Hom}(\lambda , \mu )$; it can be $\infty$, but it is
finite if $\l,\mu$ have finite index. See \cite{J1} for the
definition of index for type $II_1$ case which initiated the subject
and  \cite{PP} for  the definition of index in general. Also see
\S2.3 of \cite{KLX} for expositions. $\langle  \lambda , \mu
\rangle$ depends only on $[\lambda ]$ and $[\mu ]$. Moreover we have
if $\nu$ has finite index, then $\langle \nu \lambda , \mu \rangle =
\langle \lambda , \bar \nu \mu \rangle $, $\langle \lambda\nu , \mu
\rangle = \langle \lambda , \mu \bar \nu \rangle $ which follows
from Frobenius duality. $\mu $ is a subsector of $\lambda $ if there
is an isometry $v\in M$ such that $\mu(x)= v^* \lambda(x)v, \forall
x\in M.$ We will also use the following notation: if $\mu $ is a
subsector of $\lambda $, we will write as $\mu \prec \lambda $  or
$\lambda \succ \mu $.  A sector is said to be irreducible if it has
only one subsector.

\subsection{Local nets}
%\subsubsection{Preliminaries on conformal nets}
By an interval of the circle we mean an open connected
non-empty subset $I$ of $S^1$ such that the interior of its
complement $I'$ is not empty.
We denote by $\I$ the family of all intervals of $S^1$.

A {\it net} $\A$ of von Neumann algebras on $S^1$ is a map
\[
I\in\I\to\A(I)\subset B(\H)
\]
from $\I$ to von Neumann algebras on a fixed separable Hilbert space $\H$
that satisfies:
\begin{itemize}
\item[{\bf A.}] {\it Isotony}. If $I_{1}\subset I_{2}$ belong to
$\I$, then
\begin{equation*}
 \A(I_{1})\subset\A(I_{2}).
\end{equation*}
\end{itemize}
If $E\subset S^1$ is any region, we shall put
$\A(E)\equiv\bigvee_{E\supset I\in\I}\A(I)$ with $\A(E)=\mathbb C$
if $E$ has empty interior (the symbol $\vee$ denotes the von Neumann
algebra generated).

The net $\A$ is called {\it local} if it satisfies:
\begin{itemize}
\item[{\bf B.}] {\it Locality}. If $I_{1},I_{2}\in\I$ and $I_1\cap
I_2=\emptyset$ then
\begin{equation*}
 [\A(I_{1}),\A(I_{2})]=\{0\},
 \end{equation*}
where brackets denote the commutator.
\end{itemize}
The net $\A$ is called {\it M\"{o}bius covariant} if in addition
satisfies
the following properties {\bf C,D,E,F}:
\begin{itemize}
\item[{\bf C.}] {\it M\"{o}bius covariance}.
There exists a non-trivial strongly
continuous unitary representation $U$ of the M\"{o}bius group
$\Mob$ (isomorphic to $PSU(1,1)$) on $\H$ such that
\begin{equation*}
 U(g)\A(I) U(g)^*\ =\ \A(gI),\quad g\in \Mob,\ I\in\I.
\end{equation*}
\item[{\bf D.}] {\it Positivity of the energy}.
The generator of the one-parameter
rotation subgroup of $U$ (conformal Hamiltonian), denoted by
$L_0$ in the following,  is positive.
\item[{\bf E.}] {\it Existence of the vacuum}.  There exists a unit
$U$-invariant vector $\Omega\in\H$ (vacuum vector), and $\Omega$ is
cyclic for the von Neumann algebra $\bigvee_{I\in\I}\A(I)$.
\end{itemize}
By the Reeh-Schlieder theorem $\Omega$ is cyclic and separating for
every fixed $\A(I)$. The modular objects associated with
$(\A(I),\Omega)$ have a geometric meaning
\[
\Delta^{it}_I = U(\Lambda_I(2\pi t)),\qquad J_I = U(r_I)\ .
\]
Here $\Lambda_I$ is a canonical one-parameter subgroup of $\Mob$ and $U(r_I)$ is a
antiunitary acting geometrically on $\A$ as a reflection $r_I$ on $S^1$.

This implies {\em Haag duality}:
\[
\A(I)'=\A(I'),\quad I\in\I\ ,
\]
where $I'$ is the interior of $S^1\setminus I$.

\begin{itemize}
\item[{\bf F.}] {\it Irreducibility}. $\bigvee_{I\in\I}\A(I)=B(\H)$.
Indeed $\A$ is irreducible iff
$\Om$ is the unique $U$-invariant vector (up to scalar multiples).
Also  $\A$ is irreducible
iff the local von Neumann
algebras $\A(I)$ are factors. In this case they are either ${\mathbb C}$ or
III$_1$-factors
with separable predual
in
Connes classification of type III factors.
\end{itemize}
By a {\it conformal net} (or diffeomorphism covariant net)
$\A$ we shall mean a M\"{o}bius covariant net such that the following
holds:
\begin{itemize}
\item[{\bf G.}] {\it Conformal covariance}. There exists a projective
unitary representation $U$ of $\Diff(S^1)$ on $\H$ extending the unitary
representation of $\Mob$ such that for all $I\in\I$ we have
\begin{gather*}
 U(\phi)\A(I) U(\phi)^*\ =\ \A(\phi.I),\quad  \phi\in\Diff(S^1), \\
 U(\phi)xU(\phi)^*\ =\ x,\quad x\in\A(I),\ \phi\in\Diff(I'),
\end{gather*}
\end{itemize}
where $\Diff(S^1)$ denotes the group of smooth, positively oriented
diffeomorphism of $S^1$ and $\Diff(I)$ the subgroup of
diffeomorphisms $g$ such that $\phi(z)=z$ for all $z\in I'$. Note
that by Haag duality we have $U(\phi)\in \A(I), \forall \phi\in
\Diff (I).$ Hence the following definition makes sense:
\begin{definition}\label{virnet}
If $\A$ is a conformal net, the Virasoro subnet of $\A$, denoted by
$\Vir_\A$ is defined as follows: for each interval $I\in \I$,
$\Vir_\A(I)$ is the von Neumann algebra generated by $U(\phi)\in
\A(I), \forall \phi\in \Diff (I).$
\end{definition}
A (DHR) representation $\pi$ of $\A$ on a Hilbert space $\H$ is a
map $I\in\I\mapsto  \pi_I$ that associates to each $I$ a normal
representation of $\A(I)$ on $B(\H)$ such that
\[
\pi_{\widetilde I}\res\A(I)=\pi_I,\quad I\subset\widetilde I, \quad
I,\widetilde I\subset\I\ .
\]
$\pi$ is said to be M\"obius (resp. diffeomorphism) covariant if
there is a projective unitary representation $U_{\pi}$ of $\Mob$ (resp.
$\Diff(S^1)$) on $\H$ such that
\[
\pi_{gI}(U(g)xU(g)^*) =U_{\pi}(g)\pi_{I}(x)U_{\pi}(g)^*
\]
for all $I\in\I$, $x\in\A(I)$ and $g\in \Mob$ (resp.
$g\in\Diff(S^1)$).

By definition the irreducible conformal net is in fact an irreducible
representation of itself and we will call this representation the {\it
vacuum representation}.\par

Let $G$ be a simply connected  compact Lie group. By Th. 3.2 of
\cite{FG}, the vacuum positive energy representation of the loop
group $LG$ (cf. \cite{PS}) at level $k$ gives rise to an irreducible
conformal net denoted by {\it ${\A}_{G_k}$}. By Th. 3.3 of
\cite{FG}, every irreducible positive energy representation of the
loop group $LG$ at level $k$ gives rise to  an irreducible covariant
representation of ${\A}_{G_k}$. \par Given an interval $I$ and a
representation $\pi$ of $\A$, there is an {\em endomorphism of $\A$
localized in $I$} equivalent to $\pi$; namely $\r$ is a
representation of $\A$ on the vacuum Hilbert space $\H$, unitarily
equivalent to $\pi$, such that
$\r_{I'}=\text{id}\restriction\A(I')$. We now define  the
statistics. Given the endomorphism $\r$ of $\A$ localized in
$I\in\I$, choose an equivalent endomorphism $\r_0$ localized in an
interval $I_0\in\I$ with $\bar I_0\cap\bar I =\emptyset$ and let $u$
be a local intertwiner in $\Hom(\r,\r_0)$ , namely $u\in
\Hom(\r_{\widetilde I},\r_{0,\widetilde I})$ with $I_0$ following
clockwise $I$ inside $\widetilde I$ which is an interval containing
both $I$ and $I_0$.

The {\it statistics operator} $\epsilon (\r,\rho):= u^*\r(u) =
u^*\r_{\widetilde I}(u) $ belongs to $\Hom(\r^2_{\widetilde
I},\r^2_{\widetilde I})$. We will call $\epsilon (\r,\rho)$ the
positive or right braiding and $\widetilde\epsilon
(\r,\rho):=\epsilon (\r,\rho)^*$ the negative or left braiding. The
{\em statistics parameter} $\l_\r$ can be defined in general. In
particular, assume $\r$ to be localized in $I$ and
$\r_I\in\text{End}((\A(I))$ to be irreducible with a conditional
expectation $E: \A(I)\to \r_I(\A(I))$, then
\[
\l_\r:=E(\epsilon)
\]
depends only on the  sector of $\r$. The {\em statistical dimension}
$d_{\r}$ and the  {\it univalence} $\omega_\r$ are then defined by
\[
d_{\r} = |\lambda_\r|^{-1}\ ,\qquad \omega_\r =
\frac{\lambda_\r}{|\lambda_\r|}\ .
\]

The {\em conformal spin-statistics theorem} (cf. \cite{GL2})  shows
that
\[
\omega_\r = e^{i 2\pi L_0(\r)}\ ,
\]
where $L_0(\r)$ is the conformal Hamiltonian (the generator of the
rotation subgroup) in the representation $\r$. The right hand side
in the above equality is called the {\em univalence} of $\r$.
\par
Let $\{[\lambda], \lambda\in \L \}$ be a finite set  of all
equivalence classes of irreducible, covariant, finite-index
representations of an irreducible local conformal net $\A$. We will
denote the conjugate of $[\lambda]$ by $[{\bar \lambda}]$ and
identity sector (corresponding to the vacuum representation) by
$[1]$ if no confusion arises, and let $N_{\lambda\mu}^\nu = \langle
[\lambda][\mu], [\nu]\rangle $. Here $\langle \mu,\nu\rangle$
denotes the dimension of the space of intertwiners from $\mu$ to
$\nu$ (denoted by $\text {\rm Hom}(\mu,\nu)$).  We will denote by
$\{T_e\}$ a basis of isometries in $\text {\rm
Hom}(\nu,\lambda\mu)$. The univalence of $\lambda$ and the
statistical dimension of (cf. \S2  of \cite{GL1}) will be denoted by
$\omega_{\lambda}$ and $d{(\lambda)}$ (or $d_{\lambda})$)
respectively. The following equation is called {\it monodromy
equation} (cf.  \cite{Rehs}):\par
%Let $\phi_\lambda$ be the unique minimal left
%inverse of $\lambda$, define:
\begin{equation}\label{monodromy}
\epsilon (\mu, \lambda) \epsilon (\lambda, \mu))T_e=
\frac{\omega_\nu}{\omega_\lambda\omega_\mu} T_e
\end{equation}
where $\epsilon (\mu, \lambda)$ is the unitary braiding operator.
\par

We make the following definitions for convenience:
\begin{definition}\label{localset}
Let $\la,\mu$ be   (not necessarily irreducible) representations of
$\A$.  $H(\la,\mu):=\e(\lambda,\mu)\e(\mu,\lambda).$ We say that
$\la$ is local with $\mu$ if $H(\la,\mu)=1.$
\end{definition}
\begin{definition}\label{localsystem}
Let $\Gamma$ be a  set of DHR representations of $\A.$ If $\Gamma$
is an abelian group with multiplication given by composition and
$d_\lambda=1, \omega_\la=1, \forall \lambda\in \Gamma,$ then
$\Gamma$ is called { a local system of automorphisms}.
\end{definition}
The following Lemma  will be useful to check if a set is  a local
system of automorphims.
\begin{lemma}\label{checklocal}
(1) Assume that $[\mu]=\sum_{1\leq i\leq n} [\mu_i]$ and $\lambda,
\mu_i, i=1,...,n$ are representations of $\A.$ Then
$H(\lambda,\mu)=1$  if and only if $H(\lambda, \mu_i)=1$  for all
$1\leq i\leq n;$\par (2) If $H(\la,\mu)=1$ and $H(\la,\nu)=1$, then
$ H(\la, \mu\nu)=1$ ;\par (3) If $\la_1,...,\la_n$ generate a finite
abelian group $\Gamma$ under composition, $\omega_{\la_i}=1, 1\leq
i\leq n,$ and $H (\la_i,\la_j)=1, 1\leq i,j\leq n, $  then $\Gamma$
is a local system of automorphisms.
\end{lemma}
\proof (1) and (2) follows from \cite{Rehren} or Lemma 3.8 of
\cite{BE3}. As for (3), we prove by induction on $n.$ If $n=1,$ then
$\e(\la_1,\la_1)= \omega_{\la_1}=1$ since $\e(\la_1,\la_1)$ is a
scalar,  and it follows that $\omega_{\la_1^i}=
\e(\la_1,\la_1)^{i^2}=1, \forall i\geq 1.$
%Since
%$\omega_{\la_1^{-i}}=\omega_{\la_1^i}$ we have proved $n=1$ case.

Assume that (3) has been proved for $n-1.$ Let $\mu$ be in the
abelian group generated by $\la_1,...,\la_{n-1}.$ Since for any
integer $k$ $H(\mu, \la_n^k)=1$  by (2) and assumption, by
repeatedly applying (2) and monodromy equation, we have
$\omega_{\mu\la_n^k} = \omega_{\mu}\omega_{\la_n^k}=1$ by induction
hypotheses. It follows that (3) is proved.\endproof
\par
Next we  recall some definitions from \cite{KLM} . Recall that
${\I}$ denotes the set of intervals of $S^1$. Let $I_1, I_2\in
{\I}$. We say that $I_1, I_2$ are disjoint if $\bar I_1\cap \bar
I_2=\emptyset$, where $\bar I$ is the closure of $I$ in $S^1$. When
$I_1, I_2$ are disjoint, $I_1\cup I_2$ is called a 1-disconnected
interval in \cite{Xjw}. Denote by ${\I}_2$ the set of unions of
disjoint 2 elements in ${\I}$. Let ${\A}$ be an irreducible
M\"{o}bius covariant net . For $E=I_1\cup I_2\in{\I}_2$, let
$I_3\cup I_4$ be the interior of the complement of $I_1\cup I_2$ in
$S^1$ where $I_3, I_4$ are disjoint intervals. Let
$$
{\A}(E):= A(I_1)\vee A(I_2), \quad
\hat {\A}(E):= (A(I_3)\vee A(I_4))'.
$$ Note that ${\A}(E) \subset \hat {\A}(E)$.
Recall that a net ${\A}$ is {\it split} if ${\A}(I_1)\vee {\A}(I_2)$
is naturally isomorphic to the tensor product of von Neumann
algebras ${\A}(I_1)\otimes {\A}(I_2)$ for any disjoint intervals
$I_1, I_2\in {\I}$. ${\A}$ is {\it strongly additive} if
${\A}(I_1)\vee {\A}(I_2)= {\A}(I)$ where $I_1\cup I_2$ is obtained
by removing an interior point from $I$. \bdefin\label{rational}
\cite{{KLM}, {LX}} A M\"{o}bius covariant net ${\A}$ is said to be
completely  rational if ${\A}$ is split, strongly additive, and the
index $[\hat {\A}(E): {\A}(E)]$ is finite for some $E\in {\I}_2.$
The value of the index $[\hat {\A}(E): {\A}(E)]$ (it is independent
of $E$ by Prop. 5 of \cite{KLM}) is denoted by $\mu_{{\A}}$ and is
called the $\mu$-index of ${\A}$.
%\label{Definition 2.2}
\edefin Note that, by  results in \cite{LX}, every irreducible,
split, local conformal net with finite $\mu$-index is automatically
strongly additive. Also note that if $\A$ is completely rational,
then $\A$ has only finitely many irreducible covariant
representations by \cite{KLM}.
\par
\bdefin\label{holomorphic} A M\"{obius} net $\A$ is called
holomorphic if $\A$ is completely rational and $\mu_\A=1,$ i.e.,
$\A$ has only one irreducible representation which is the vacuum
representation. \edefin

Let $\B$ be a  M\"{o}bius (resp. conformal) net. $\B$ is called a
M\"{o}bius (resp. conformal) extension of $\A$ if there is  a map
\[
I\in\I\to\A(I)\subset \B(I)
\]
that associates to each interval $I\in \I$ a von Neumann subalgebra $\A(I)$
of $\B(I)$, which is isotonic
\[
\A(I_1)\subset \A(I_2), I_1\subset I_2,
\]
and   M\"{o}bius (resp. diffeomorphism) covariant with respect to
the representation $U$, namely
\[
U(g) \A(I) U(g)^*= \A(g.I)
\] for all $g\in \Mob$ (resp. $g\in \Diff(S^1)$) and $I\in \I$. $\A$ will be called a  M\"{o}bius (resp. conformal)
subnet of $\B.$ Note that by Lemma 13 of \cite{L1} for each $I\in
\I$ there exists a conditional expectation $E_I: \B(I)\rightarrow
\A(I)$ such that $E$ preserves the vector state given by the vacuum
of $\B$.
\begin{definition}\label{ext}
Let $\A$ be a  M\"{o}bius covariant net. A  M\"{o}bius covariant net
$\B$ on a Hilbert space $\H$ is an extension of $\A$ if there is a
DHR representation  $\pi$ of $\A$ on $\H$ such that $\pi(\A)\subset
\B$ is a   M\"{o}bius subnet. The extension is irreducible if
$\pi(\A(I))'\cap \B(I) = {\mathbb C} $ for some (and hence all)
interval $I$, and is of finite index if $\pi(\A(I))\subset \B(I)$
has finite index for some (and hence all) interval $I$. The index
will be called the index of the inclusion $\pi(\A)\subset \B$ and
will be denoted by $[\B:\A].$  If $\pi$ as representation of
 $\A$ decomposes as $[\pi]= \sum_\lambda m_\lambda[\lambda]$ where
$m_\lambda$  are non-negative  integers and $\lambda$ are irreducible
DHR representations of $\A$, we say that
$[\pi]= \sum_\lambda m_\lambda[\lambda]$ is the spectrum of the
extension.  For simplicity we will write  $\pi(\A)\subset \B$ simply as
$\A\subset \B$.
\end{definition}
\begin{lemma}\label{indexab}
If  $\A$ is completely rational, and a   M\"{o}bius covariant net
$\B$ is an irreducible extension of $\A$. Then $\A\subset\B$ has
finite index , $\B$ is completely rational and
$$\mu_\A= \mu_\B [\B:\A]^2.$$
\end{lemma}
\proof
 $\A\subset\B$ has finite
index follows from Prop. 2.3 of \cite{KL}, and the rest follows from
Prop. 24 of \cite{KLM}.
\endproof
\begin{lemma}\label{extconformal}
If  $\A$ is a conformal net, and a   M\"{o}bius covariant net $\B$
is an  extension of $\A$ with index $[\B:\A]< \infty.$ Then $\B$ is
a conformal net.\end{lemma} \proof Denote by $\pi$ the vacuum
representation of $\B.$ Denote by $\bold G$  the universal cover of
$\Mob$. By definition $g\in {\bold G}\rightarrow U_\pi(g)$ is a
representation of ${\bold G}$ which implements the M\"{o}bius
covariance of $\pi\res \A.$ On the other hand by \S2 of \cite{AFK}
there is a representation of $g\in {\bold G}\rightarrow V_\pi(g)$
which implements the M\"{o}bius covariance of $\pi\res \A,$ and
$V_\pi(g)\in \bigvee_{I\in \I} \pi(\Vir_\A(I)),$ where $\Vir_\A$ is
defined in definition \ref{virnet}.  Since by assumption $\pi\res
\A$ has finite index, by Prop. 2.2 of \cite{GL1} we have
$U_\pi(g)=V_\pi(g), \forall g\in {\bold G}.$ Hence $\Vir_\A\subset
\B$ verifies the condition in definition 3.1 of \cite{Ca}, and by
Prop. 3.7 of \cite{Ca} the lemma is proved.
\endproof

The following is Th. 4.9 of \cite{LR} (cf. \S2.4 of \cite{KL}) which
is also used in \S4.2 of  \cite{KL}:
\begin{proposition}\label{qlocal}
Let $\A$ be a  M\"{o}bius covariant net, $\r$ a DHR representation
of $\A$ localized on a fixed $I_0$ with finite statistics, which contains
$\id$ with multiplicity one, i.e., there is (unique up to a phase)
isometry $w\in\Hom(\id,\r).$ Then there is a  M\"{o}bius
covariant net $\B$ which is an irreducible extension of $\A$ if and only
if there is an isometry $w_1\in\Hom(\r,\r^2)$
which solves the following equations:
\begin{align}
w_1^* w & = w_1^* \r(w) \in {\mathbb R_+} \label{a}\\
w_1w_1& = \r(w_1) w_1 \label{b} \\
\epsilon(\r,\r) w_1 & = w_1 \label{c}
\end{align}
\end{proposition}
\begin{remark}\label{outer}
Let $\A\subset \B$ be as in Prop.\ref{qlocal}. If $U$ is an unitary
on the vacuum representation space of $\A$ such that $\Ad_U \A(I)=
\A(I), \forall I,$ then it is easy to check that $(\Ad_U
\rho\Ad_U^*, \Ad_U (w_1), \Ad_U (w))$ verifies the equations in
Prop. \ref{qlocal}, and  determines a  M\"{o}bius covariant net
$\Ad_U(\B)$. The spectrum of $\A\subset\Ad_U(\B)$ (cf. definition
\ref{ext}) is $\Ad_U \rho\Ad_U^*$ which may be different from
$\rho$, but $\Ad_U(\B)$ is isomorphic to $\B$ by definition.
\end{remark}
\subsection{Induction}\label{ind}
Let $\B$ be a M\"obius covariant net and $\A$ a subnet. We assume
that $\A$ is strongly additive and $\A\subset \B$ has finite index.
Fix an interval $I_0\in\I$ and  canonical endomorphism (cf.
\cite{LR}) $\gamma$ associated with $\A(I_0)\subset\B(I_0)$. Then we
can choose for each $I\subset\I$ with $I\supset I_0$ a canonical
endomorphism $\gamma_{I}$ of $\B(I)$ into $\A(I)$ in such a way that
$\gamma_{I}\res\B(I_0)=\gamma_{I_0}$ and $\rho_{I_1}$ is the
identity on $\A(I_1)$ if $I_1\in\I_0$ is disjoint from $I_0$, where
$\rho_{I}\equiv\gamma_{I}\res\A(I)$. Given a DHR endomorphism $\la$
of $\A$ localized in $I_0$, the inductions $\a_{\la},\a_{\la}^{-}$
of $\la$ are the endomorphisms of $\B(I_0)$ given by
\[
\a_{\la}\equiv \gamma^{-1}\cdot\Ad\e(\la,\rho)\cdot\la\cdot\gamma \
, \a_{\la}^{-}\equiv
\gamma^{-1}\cdot\Ad\tilde{\e}(\la,\rho)\cdot\la\cdot\gamma
\]
where $\e$ (resp. $\tilde\e$) denotes the right braiding (resp. left
braiding) (cf. Cor. 3.2 of \cite{BE1}). In \cite{Xb} a slightly
different endomorphism was introduced and the relation between the
two was given in \S2.1 of \cite{X3m}.

Note that $\Hom( \a_\lambda,\a_\mu)=:\{ x\in \B(I_0) | x
\a_\lambda(y)= \a_\mu(y)x, \forall y\in \B(I_0)\} $ and $\Hom(
\lambda,\mu)= :\{ x\in \A(I_0) | x  \lambda(y)= \mu(y)x, \forall
y\in \A(I_0)\} .$

The following is Lemma 3.6 of \cite{BE3} and Lemma 3.5 of
\cite{BE1}:
\begin{lemma}\label{aa'}
$$\Hom (\a_\lambda^{}, \a_\mu^{-}) = \{ T \in \B(I_0) | \gamma(T) \in \Hom (\rho \lambda,
\rho\mu) | \e(\mu,\rho)\e(\rho,\mu)\gamma(T)=\gamma(T)\}.
$$
\end{lemma}
As a consequence of Lemma \ref{aa'} we have the following Prop. 3.23
of \cite{BE1} ( Also cf. the proof of Lemma 3.2 of \cite{Xb}):
\begin{lemma}\label{a=a'}
$[\a_\lambda^{}]=[\a_\lambda^{-}]$ iff $\e(\lambda, \rho) \e(
\rho,\la)=1$ .
\end{lemma}

The following follows from Lemma 3.4 and Th. 3.3 of \cite{Xb} (also
cf. \cite{BE1}) :
\begin{lemma}\label{3.3}
(1): $[\la]\rightarrow [\a_\la], [\la]\rightarrow [\a_\la^{-}]$ are
ring homomorphisms;\par (2) $\lan \a_\la,\a_\mu\ran = \lan \la \rho,
\mu\ran.$
\end{lemma}

\subsection{Local simple current extensions}\label{simpleextension}

\begin{proposition}\label{simple}
(1) Assume that $\B$ is a M\"{o}bius extension of $\A$ of finite
index with spectrum $[\pi]=\sum_{\lambda \in \exp} m_\l[\l].$  Let
$\Gamma:=\{ \l| \l\in \exp\}.$ Assume that $d_\l=1, \forall \l\in
\exp.$ Then $\Gamma$ is a local system of automorphisms;\par (2) If
$\Gamma$ is a finite local system of automorphisms of $\A$, then
there is a M\"{o}bius extension $\B$ of $\A$ with spectrum
$[\pi]=\sum_{\lambda \in \Gamma} [\l].$
\end{proposition}
\proof Ad (1): By assumption we have $\a_\la^{}\succ 1, \forall
\la\in \Gamma.$ By Lemma 3.10 of \cite{BE3} $\omega_\la=1.$ Since
$d_\l=d_{\a_\la}=1,$ it follows that $[\a_\la]=[\a_\la^{-}]=[1].$
Note that if $\la\in \Gamma$ iff $[\a_\la]=[1]$ and it follows that
$\Gamma$ is an abelian group with multiplication given by
composition. By Lemma \ref{a=a'} and Lemma \ref{checklocal} (1) is
proved.\par (2) It follows from Prop. 5.5 of \cite{Rehren} (also. cf
Th. 5.2 of \cite{DHR}) that there is a M\"{o}bius extension $\B$ of
$\A$ with spectrum $[\pi]=\sum_{\lambda \in \Gamma} [\l].$
\endproof

\begin{remark}
(1) We will use the notation $\B=\A\ltimes \Gamma$ for the extension
in Prop. \ref{simple}.

(2) One can extend the above theorem to a case when $\B$ is not
local but verifies twisted locality. Such extensions have been used
for example in \cite{KLPR}.
\end{remark}

\subsection{Mirror extensions}\label{mirrorextension}
%\subsection{Coset construction}
In this section we recall the mirror construction as given in \S3 of
\cite{Xm}.  Let $\B$ be a completely rational net and $\A\subset \B$
be a subnet which is also completely rational.
\begin{definition}\label{coset}
Define a subnet $\widetilde\A\subset \B$ by $\widetilde\A(I):=
\A(I)'\cap \B(I), \forall I\in \I.$
\end{definition}
We note that since $\A$ is completely rational, it is strongly
additive and so we have $\widetilde\A(I)= (\vee_{J\in \I}\A(J))'\cap
\B(I), \forall I\in \I.$ The following lemma then follows directly
from the definition:
\begin{lemma}\label{cosetnet}
The restriction of $\widetilde\A$ on the Hilbert space $ \overline{
\vee_I\widetilde\A(I)\Omega}$ is an irreducible
  M\"{o}bius  covariant net.
\end{lemma}
The net $\widetilde\A$ as in Lemma \ref{cosetnet} will be called the
{\it coset} of  $\A\subset \B$. See \cite{Xcos} for a class of
cosets from Loop groups. \par The following definition generalizes
the definition in \S3 of \cite{Xcos}:
\begin{definition}\label{cofinite}
 $\A\subset \B$ is called cofinite if the inclusion
$\widetilde\A(I)\vee \A(I) \subset \B(I)$ has finite index for some
interval $I$.
\end{definition}
The following is Prop. 3.4 of \cite{Xm}:
\begin{proposition}\label{rationalc}
Let $\B$ be completely rational, and let $\A\subset\B$ be a
M\"{o}bius subnet which is also completely rational. Then $\A\subset
\B$ is  cofinite if and only if $\tilde\A$ is completely rational.
\end{proposition}
Let  $\B$ be completely rational, and let $\A\subset\B$ be a
M\"{o}bius subnet which is also completely rational. Assume that
$\A\subset \B$ is cofinite. We will use $\sigma_i,\sigma_j,...$
(resp. $\lambda, \mu...$) to label irreducible DHR representations
of $\B$ (resp. $\A$) localized on a fixed interval $I_0$. Since $
\widetilde\A$ is completely rational by Prop. \ref{rationalc},
$\widetilde\A\otimes \A$ is completely rational, and so every
irreducible DHR representation $\sigma_i$ of $\B$, when restricting
to $\widetilde\A\otimes \A$, decomposes as direct sum of
representations of  $ \widetilde\A\otimes \A$ of the form
$(i,\lambda)\otimes \lambda$ by Lemma 27 of \cite{KLM}. Here
$(i,\lambda)$ is a DHR representation of $\widetilde \A$ which may
not be irreducible and we use the  tensor notation
$(i,\lambda)\otimes \lambda$ to represent a DHR representation of $
\widetilde\A\otimes \A$ which is localized on $I_0$ and defined by
$$
(i,\lambda)\otimes \lambda (x_1\otimes x_2)= (i,\lambda)(x_1)\otimes
 \lambda(x_2), \forall x_1\otimes x_2\in \widetilde\A(I_0) \otimes \A(I_0).
$$
We will also identify $\widetilde \A$ and $\A$ as subnets of $
\widetilde\A\otimes \A$ in the natural way. We note that when no
confusion arise, we will use $1$ to denote the vacuum representation
of a net.

\begin{definition}\label{normal}
A  M\"{o}bius subnet $\A\subset\B$ is normal if $\widetilde
\A(I)'\cap \B(I)= \A$ for some I.
\end{definition}
The following is implied by Lemma 3.4 of \cite{Reh1} (also cf. Page
797 of \cite{Xcos2}):
\begin{lemma}\label{normal1}
Let  $\B$ be completely rational, and let $\A\subset\B$ be a
M\"{o}bius subnet which is also completely rational. Assume that
$\A\subset \B$ is cofinite. Then the following conditions are
equivalent: \par (1) $\A\subset\B$ is normal; \par (2) $(1,1)$ is
the vacuum representation of $\widetilde \A$ and $(1,\lambda) $
contains $(1,1)$ if and only if $\lambda=1$. \par
\end{lemma}
The following is part of Proposition 3.7 of \cite{Xm}:
\begin{proposition}\label{mirror}
Let  $\B$ be completely rational, and let $\A\subset\B$ be a
M\"{o}bius subnet which is also completely rational. Assume that
$\A\subset \B$ is cofinite and normal. Then:\par (1) Let $\gamma$ be
the restriction of the vacuum representation of $\B$ to
$\widetilde\A\otimes \A$. Then $[\gamma]= \sum_{\lambda\in \exp}
[(1,\lambda)\otimes \lambda]$ where each $(1,\lambda)$ is
irreducible;\par (2) Let $\lambda\in \exp$ be as in (1), then
$[\alpha_{(1,\lambda)\otimes 1}] = [\alpha_{1\otimes \bar\lambda}]$,
and $[\lambda]\rightarrow [\alpha_{1\otimes \lambda}]$ is a ring
isomorphism where the $\a$-induction is with respect to
$\widetilde\A\otimes \A \subset \B$ as in subsection \ref{ind};
Moreover the set $\exp$ is closed under fusion;
\par
(3) Let $[\rho]= \sum_{\lambda\in \exp} m_\lambda[\lambda]$ where
$m_\lambda = m_{\bar\lambda}\geq 0, \forall \lambda$, and $
[(1,\rho)]= \sum_{\lambda\in \exp} m_\lambda[(1,\lambda)]$. Then there
exists an unitary element $T_\rho \in
\Hom(\alpha_{(1,\rho)\otimes 1},
\alpha_{1\otimes \rho})$ such that
$$
\epsilon((1,\rho),(1,\rho)) T_\rho^* \alpha_{1\otimes
\rho}(T_\rho^*) =T_\rho^* \alpha_{1\otimes \rho}(T_\rho^*)
\widetilde\epsilon(\rho,\rho) ;$$ \par (4) Let $\rho$, $(1,\rho)$ be
as in (3). Then
\begin{align*}
\Hom(\rho^n, \rho^m) & = \Hom(\alpha_{1\otimes \rho^n},\alpha_{1\otimes \rho^m}),  \\
\Hom((1,\rho)^n, (1,\rho)^m) & = \Hom(\alpha_{(1,\rho)^n\otimes 1},
\alpha_{(1,\rho)^m\otimes 1}), \forall n,m\in {\mathbb N};
\end{align*}
\end{proposition}
Denote by $\Delta_0:=\{ \la| [\la]=\sum_i{[\la_i]}, \la_i\in exp\}.$
Assume  $\mu_i \in \Delta_0, i=1,...,n.$ For each $[\mu_i]=\sum_j
m_{ij}[\la_j],$ choose DHR representations $M(\mu_i)$ of
$\widetilde\A$ such that $[M(\mu_i)]=\sum_j m_{ij}[(1,\la_j)].$ Let
$T_i\in \Hom(\a_{ m(\mu_i)\otimes 1},\a_{1\otimes \mu_i})$ be an
unitary element (not necessarily unique up to phase when $\mu_i$ is
not irreducible) as given is (3) of Prop. \ref{mirror}. Define
$$
T_{i_1i_2...i_k}:=\a_{\mu_1...\mu_{k-1}\otimes 1}(T_{i_k})...
\a_{\mu_1...\mu_{2}\otimes 1}(T_{i_3}) \a_{\mu_1\otimes
1}(T_{i_2})T_{i_1}\in \Hom(\a_{ M(\mu_1)...M(\mu_k)\otimes 1},
\a_{1\otimes \mu_1...\mu_k})
$$
For each $S\in \Hom(\mu_{i_1}...\mu_{i_k}, \mu_{j_1}...\mu_{j_m})$
we define $M(S):= T^*_{j_1...j_m} S T_{i_1...i_k}.$
\begin{lemma}\label{M}
Assume that  $S_1, T\in \Hom(\la,\mu), S_2\in \Hom(\nu,\la)$ where
$\la,\mu$ are products of elements from $\{ \mu_1,...,\mu_n\}.$ If
$\nu=\mu_{i_1}...\mu_{i_k}$ we denote $M(\mu_{i_1})... M(\mu_{i_k})$
by $M(\nu).$ Then:
$$
M(S_1S_2)=M(S_1)M(S_2), M(\nu (T))= M(\nu)(M(T)),
M(\e(\mu_i,\mu_j))= \tilde{\e}(M(\mu_i),M(\mu_j)).
$$
\end{lemma}
\proof The first two identities follow directly from definitions.
The third follows from (3) of Prop. 2.3.1 of \cite{X3m}, as (3) of
Prop. \ref{mirror}.
\endproof
The following is Th. 3.8 of \cite{Xm}:
\begin{theorem}\label{mainmirror} Let  $\B$ be completely rational, and let
$\A\subset\B$ be a  M\"{o}bius subnet which is also completely
rational. Assume that  $\A\subset \B$ is cofinite and normal, and
let $\exp$ be as in (1) of Prop.\ref{mirror}.  Assume that
$\A\subset\C$ is an irreducible   M\"{o}bius extension of $\A$ with
spectrum  $[\rho]=\sum_{\lambda\in \exp} m_\lambda [\lambda],
m_\lambda\geq 0.$ Then there is  an irreducible   M\"{o}bius
extension $\widetilde \C$ of $\widetilde\A$ with spectrum
$[(1,\rho)]=\sum_{\lambda\in \exp} m_\lambda [(1,\lambda)]$.
Moreover $\widetilde \C$ is completely rational.
\end{theorem}
\begin{remark}\label{mc}
Due to (5) of Prop. 3.7 of \cite{Xm}, the extension $\widetilde
\A\subset\widetilde \C$ as given in Th. \ref{main} will be called
the mirror or the conjugate of $\A\subset \C$.
\end{remark}

 By Lemma \ref{indexab} and Th. \ref{mainmirror} we have:
\begin{corollary}\label{mirrorindex}
Let ${\widetilde \C} $ be the mirror extension as given in Th.
\ref{mainmirror}. Then
$\frac{\mu_{\widetilde{\C}}}{\mu_{\widetilde{\A}}}=
\frac{\mu_{{\C}}}{\mu_{\A}}.$
\end{corollary}
The mirror extension $\widetilde \A\subset\widetilde{C}$ is
constructed as follows: let $(\rho,w,w_1)$ be associated with
extension $\A\subset \C$ as given in Prop. \ref{qlocal}. Then the
extension $\widetilde \A\subset\widetilde{C}$ is given by $(M(\rho),
M(w), M(w_1))$ where the map $M$ is defined before Lemma \ref{M}.
Let $\mu,\nu\in \Delta_0.$ Consider now inductions with respect to
$\A\subset \C$ and $\widetilde \A\subset\widetilde{C}.$

\begin{proposition}\label{mirrora}
Assume that $\mu,\nu\in \Delta_0$, $M(\rho)=(1,\rho),
M(\mu)=(1,\mu), M(\nu)=(1,\nu).$ Then
$$
\lan \a_{\mu}, \a_{\nu}\ran =\lan \a_{M(\mu)}, \a_{M(\nu)}\ran,
 \lan \a_{\mu}, \a_{\nu}^{-}\ran =\lan
\a_{M(\mu)}, \a_{M(\nu)}^{-}\ran.
$$
\end{proposition}
\proof By Lemma \ref{3.3}  we have
$$
\lan \a_{\mu}, \a_{\nu}\ran = \lan\rho \mu,\nu\ran, \lan
\a_{M(\mu)}, \a_{M(\nu)}\ran = \lan M(\rho) M(\mu),M(\nu)\ran,
$$
By Prop. \ref{mirror} we have proved the first equality. By Lemma
\ref{aa'} we have
$$
\Hom (\a_{\mu}, \a_{\nu}^{-}) = \{ T\in \C(I_0)| \gamma (T)\in \Hom
(\rho\mu, \rho\nu), \e(\nu,\rho)\e(\rho,\nu)\gamma (T)= \gamma(T)\}.
$$
By \cite{LR}, $\gamma (\C(I_0))= \{ x\in \A(I_0)| x= w_1^* \rho(x)
w_1 \}.$ It follows that $\lan \a_{\mu}^{+}, \a_{\nu}^{-}\ran$ is
equal to the dimension of the following vector space
$$\{ T'\in \A(I_0) |
T'\in \Hom (\rho\mu, \rho\nu), \e(\nu,\rho)\e(\rho,\nu) T'= T', T'=
w_1^* \rho(T') w_1 \}.
$$
Now apply the map $M$ to the above vector space  and use Lemma
\ref{M} we have that $\lan \a_{\mu}^{+}, \a_{\nu}^{-}\ran$ is equal
to the dimension of the following vector space
\begin{align*}
\{ T'\in \widetilde{\A(I_0)} | T'\in \Hom (M(\rho)M(\mu),
M(\rho)M(\nu)), \tilde{\e}(M(\nu),M(\rho))\tilde{\e}(M(\rho),M(\nu))
T'= T',\\
T'= M(w_1)^* M(\rho)(T') M(w_1) \}.
\end{align*}

Since $\tilde{\e}(M(\nu),M(\rho))\tilde{\e}(M(\rho),M(\nu))=
(\e(M(\nu),M(\rho))\e(M(\rho),M(\nu)))^*$, we conclude that $\lan
\a_{\mu}^{+}, \a_{\nu}^{-}\ran$ is equal to
 the dimension of the following vector space
\begin{align*}
\{ T'\in \widetilde{\A(I_0)} | T'\in \Hom (M(\rho)M(\mu),
M(\rho)M(\nu)), {\e(M(\nu),M(\rho))}{\e(M(\rho),M(\nu))} T'= T',
\\
T'= M(w_1)^* M(\rho)(T') M(w_1) \}
\end{align*}
which is equal to $\lan \a_{ M(\mu)},\a_{M(\nu)}^-\ran$ by Lemma
\ref{aa'}.
\endproof
\subsection{A series of normal extensions}\label{lr2}
 Let $G= SU(n)$. We denote $LG$ the group of
smooth maps $f: S^1 \mapsto G$ under pointwise multiplication. The
diffeomorphism group of the circle $\text{\rm Diff} S^1 $ is
naturally a subgroup of $\text{\rm Aut}(LG)$ with the action given
by reparametrization. In particular the group of rotations
$\text{\rm Rot}S^1 \simeq U(1)$ acts on $LG$. We will be interested
in the projective unitary representation $\pi : LG \rightarrow U(H)$
that are both irreducible and have positive energy. This means that
$\pi $ should extend to $LG\ltimes \text{\rm Rot}\ S^1$ so that
$H=\oplus _{n\geq 0} H(n)$, where the $H(n)$ are the eigenspace for
the action of $\text{\rm Rot}S^1$, i.e., $r_\theta \xi = \exp(i n
\theta)$ for $\theta \in H(n)$ and $\text{\rm dim}\ H(n) < \infty $
with $H(0) \neq 0$. It follows from \cite{PS} that for fixed level
$k$ which is a positive integer, there are only finite number of
such irreducible representations indexed by the finite set
$$
 P_{++}^{k}
= \bigg \{ \lambda \in P \mid \lambda = \sum _{i=1, \cdots , n-1}
\lambda _i \Lambda _i , \lambda _i \geq 0\, , \sum _{i=1, \cdots ,
n-1} \lambda _i \leq k \bigg \}
$$
where $P$ is the weight lattice of $SU(n)$ and $\Lambda _i$ are the
fundamental weights. We will write $\la=(\la_1,...,\la_{n-1}),
\la_0= k-\sum_{1\leq i\leq n-1} \la_i$ and refer to
$\la_0,...,\la_{n-1}$ as components of $\la.$

We will use $k\Lambda_0$ or simply $1$  to denote the trivial
representation of $SU(n)$.

For $\lambda , \mu , \nu \in P_{++}^{k}$, define $N_{\lambda
\mu}^\nu  = \sum _{\delta \in P_{++}^{k} }S_\lambda ^{(\delta)}
S_\mu ^{(\delta)} S_\nu ^{(\delta*)}/S_{\Lambda_0}^{(\delta})$ where
$S_\lambda ^{(\delta)}$ is given by the Kac-Peterson formula:
$$
S_\lambda ^{(\delta)} = c \sum _{w\in S_n} \varepsilon _w \exp
(iw(\delta) \cdot \lambda 2 \pi /n)
$$
where $\varepsilon _w = \text{\rm det}(w)$ and $c$ is a
normalization constant fixed by the requirement that
$S_\mu^{(\delta)}$ is an orthonormal system. It is shown in
\cite{Kac2} P. 288 that $N_{\lambda \mu}^\nu $ are non-negative
integers. Moreover, define $ Gr(C_k)$ to be the ring whose basis are
elements of $ P_{++}^{k}$ with structure constants $N_{\lambda
\mu}^\nu $.
  The natural involution $*$ on $ P_{++}^{k}$ is
defined by $\lambda \mapsto \lambda ^* =$ the conjugate of $\lambda
$ as representation of $SU(n)$.\par

We shall also denote $S_{\Lambda _0}^{(\Lambda)}$ by $S_1^{(\Lambda
)}$. Define $d_\lambda = \frac {S_1^{(\lambda )}}{S_1^{(\Lambda
_0)}}$. We shall call $(S_\nu ^{(\delta )})$ the $S$-matrix of
$LSU(n)$ at level $k$. \par
 We shall encounter the $\Bbb Z_n$
group of automorphisms of this set of weights, generated by
$$
J : \lambda = (\lambda_1, \lambda_2, \cdots , \lambda_{n-1})
\rightarrow J(\lambda) = ( k -1- \lambda_1 -\cdots \lambda_{n-1},
\lambda_1, \cdots , \lambda_{n-2}).
$$
We will identity $J$ with $k\La_1$ in the following. Define
$\col(\lambda) = \Sigma_i (\lambda_i - 1) i $. The central element $
\exp \frac{2\pi i}{n}$ of $SU(n)$ acts on representation of $SU(n)$
labeled by $\lambda$ as $\exp( \frac{2\pi i \col(\lambda)}{n})$.
modulo $n$ $\col(\la)$ will be called the color of $\la.$

The irreducible positive energy representations of $ L SU(n)$ at
level $k$ give rise to an irreducible conformal net $\A_{SU(n)_k}$
(cf. \cite{KLX}) and its covariant representations.  $\A_{SU(n)_k}$
is completely rational (cf. \cite{W} and \cite{Xjw}), and
$\mu_{\A_{SU(n)_k}}= \frac{1}{(S_{1}^1)^2}$ by \cite{Xjw}.  We will
use $\la=(\la_1,...\la_{n-1})$ to denote irreducible representations
of $\A$ and also the corresponding endomorphism of $M=\A(I).$

All the sectors $[\lambda]$ with $\lambda$ irreducible generate the
fusion ring of $\A.$
\par For $\lambda$ irreducible, the univalence $\omega_\lambda$ is given
by an explicit formula (cf. 9.4 of [PS]). Let us first define
$h_\lambda = \frac {c_2(\lambda)}{k+n}$ where $c_2(\lambda)$ is the
value of Casimir operator on representation of $SU(n)$ labeled by
dominant weight $\lambda$.
 $h_\lambda$ is usually called the conformal dimension. Then
we have: $\omega_\lambda = exp({2\pi i} h_\lambda)$. The conformal
dimension of $\lambda=(\la_1,...,\la_{n-1})$ is given by
\begin{equation}\label{cdim} h_\lambda= \frac{1}{2n(k+n)}\sum_{1\leq
i\leq n-1} i(n-i) \la_i^2 + \frac{1}{n(k+n)}\sum_{1\leq j\leq i\leq
n-1}j (n-i)\la_j\la_i + \frac{1}{2(k+n)}\sum_{1\leq j\leq n-1}
j(n-j) \la_j \end{equation}

\par Let $G \subset H$ be inclusions of compact
simple Lie groups.  $LG \subset LH$ is called a conformal inclusion
if the level 1 projective positive energy representations of $LH$
decompose as a finite number of irreducible projective
representations of $LG$.   $LG \subset LH$ is called a maximal
conformal inclusion if there is no proper subgroup $G'$ of $H$
containing $G$ such that   $LG \subset LG'$ is also  a conformal
inclusion. A list of maximal conformal inclusions can be found in
\cite{GNO}.  \par Let $H^0$ be the vacuum representation of $LH$,
i.e., the
 representation of $LH$ associated with the trivial  representation of $H$.
Then $H^0$ decomposes as a direct sum of  irreducible projective
representation of $LG$ at level $K$. $K$  is called the Dynkin index
of the conformal inclusion.

%Assume $H^0 = \oplus_{\lambda \in P_0} m_\lambda
%H_\lambda$ where $P_0 \subset C_K$ is finite and $m_\lambda$ is the
%multiplicity of $H_\lambda$ in $H^0$ .
We shall  write the conformal inclusion as $G_K\subset H_1$. Note
that it follows from the definition that $\A_{H_1}$ is an extension
of $\A_{G_K}$.
We will be interested in the following  conformal inclusion:
$$
L(SU(m)_n \times SU(n)_m)  \subset \ L \ SU(nm)
$$
In the classification of conformal inclusions in [GNO], the above
conformal inclusion corresponds to the Grassmanian
$SU(m+n)/SU(n)\times SU(m)\times U(1)$.\par Let $\Lambda_0$ be the
vacuum representation of $LSU(nm)$ on Hilbert space $H^0$.  The
decomposition of $\Lambda_0$ under $L(SU(m) \times SU(n))$ is known,
see, e.g. \cite{ABI}.  To describe such a decomposition, let us
prepare some notation.  We  use $\dot S$ to denote the $S$-matrices
of $SU(m)$, and $\ddot S$ to denote the $S$-matrices of $SU(n)$. The
level $n$ (resp. $m$) weight of $LSU(m)$ (resp. $LSU(n)$) will be
denoted by $\dot \lambda$ (resp. $\ddot \lambda$). \par We start by
describing $\dot P_+^n$ (resp. $\ddot P_+^m$), i.e. the highest
weights of level $n$ of $LSU(m)$ (resp. level $m$ of $LSU(n)$).

$\dot P_+^n$ is the set of weights
$$
\dot \lambda = \widetilde k_0 \dot \Lambda_0 + \widetilde k_1 \dot
\Lambda_1 + \cdots + \widetilde k_{m-1} \dot \Lambda_{m-1}
$$
where $\widetilde k_i$ are non-negative integers such that
$$
\sum_{i=0}^{m-1} \widetilde k_i = n
$$
and $\dot \Lambda_i = \dot \Lambda_0 + \dot \omega_i$, $1 \leq i
\leq m-1$, where $\dot \omega_i$ are the fundamental weights of
$SU(m)$.

Instead of $\dot \lambda$ it will be more convenient to use
$$
\dot \lambda + \dot \rho = \sum_{i=0}^{m-1} k_i \dot \Lambda_i
$$
with $k_i = \widetilde k_i + 1$ and $\overset m-1 \to{\underset i=0
\to \sum} k_i = m + n$.  Due to the cyclic symmetry of the extended
Dynkin diagram of $SU(m)$, the group $\Bbb Z_m$ acts on $\dot P_+^n$
by
$$
\dot \Lambda_i \rightarrow \dot \Lambda_{(i+ \dot \mu)\mod m}, \quad
\dot \mu \in \Bbb Z_m.
$$
Let $\Omega_{m,n} = \dot P_+^n / \Bbb Z_m$.  Then there is a natural
bijection between $\Omega_{m,n}$ and $\Omega_{n,m}$ (see \S2 of
\cite{ABI}).  The idea is to draw a circle and divide it into $m+n$
arcs of equal length.  To each partition $\sum_{0\leq i\leq m-1} k_i
= m+n$ there corresponds a "slicing of the pie" into $m$ successive
parts with angles $2\pi k_i/(m+n)$, drawn with solid lines. We
choose this slicing to be clockwise.  The complementary slicing in
broken lines (The lines which are not solid) defines a partition of
$m+n$ into $n$ successive parts, $\sum_{0\leq i\leq n-1} l_i = m+n$.
We choose the later slicing to be counterclockwise, and it is easy
to see that such a slicing corresponds uniquely to an element of
$\Omega_{n,m}$.
\par

We  parameterize the bijection by a map
$$
\beta : \dot P_+^n \rightarrow \ddot P_+^m
$$
as follows.  Set
$$
r_j = \sum^m_{i=j} k_i, \quad 1 \leq j \leq m
$$
where $k_m \equiv k_0$.  The sequence $(r_1, \ldots , r_m)$ is
decreasing, $m + n = r_1 > r_2 > \cdots > r_m \geq 1$.  Take the
complementary sequence $(\bar r_1, \bar r_2, \ldots , \bar r_n)$ in
$\{ 1, 2, \ldots , m+n \}$ with $\bar r_1 > \bar r_2 > \cdots > \bar
r_n$.  Put
$$
S_j = m + n + \bar r_n - \bar r_{n-j+1}, \quad 1 \leq j \leq n.
$$
Then $m + n = s_1 > s_2 > \cdots > s_n \geq 1$.  The map $\beta$ is
defined by
$$
(r_1, \ldots , r_m) \rightarrow (s_1, \ldots , s_n).
$$
The following lemma summarizes what we will use:
\begin{lemma}\label{lr}
%\rm{(1)}
(1) Let $\dot Q$ be the root lattice of $SU(m), \ \dot \Lambda_i, \
0 \leq i \leq m-1$ its fundamental weights and $\dot Q_i = (\dot Q +
\dot \Lambda_i) \cap \dot P_+^n$. Let $\Lambda \in {\Bbb Z_{mn}}$
denote a level 1 highest weight of $SU(mn)$ and $\dot \lambda \in
\dot Q_{\Lambda \text{\rm mod} m}$.
 Then there exists a unique $\ddot \lambda
\in \ddot P_+^m$ with $\ddot \lambda = \mu \beta(\dot \lambda)$ for
some unique $\mu \in \Bbb Z_n$ such that $H_{\dot \lambda} \otimes
H_{\ddot \lambda}$ appears once and only once in $H^\Lambda$. The
map $\dot \lambda \rightarrow \ddot \lambda = \mu \beta(\dot
\lambda)$ is one-to-one. Moreover, $H^\Lambda$, as representations
of $L(SU(m) \times SU(n))$, is a direct sum of all such $H_{\dot
\lambda} \otimes H_{\ddot \lambda}$;\par (2) $\mu_{\A_{SU(n)_m}}=
\frac{n}{m} \mu_{\A_{SU(m)_n}};$ \par (3) The subnets
$\A_{{SU}(n)_m}\subset \A_{{SU}(nm)_1}$ are normal and cofinite. the
set $\exp$ as in (1) Prop. \ref{mirror} is the elements of
$P_{++}^{n+m}$ which belong to the root lattice of ${SU}(n)$.
\end{lemma}
\proof (1) is Th. 1 of \cite{ABI}. (2) follows from Th. 4.1 of
\cite{Xjw}. (3) is Lemma 4.1 of \cite{Xm}.
\endproof
\section{Schellekens's modular invariants and their realizations by conformal nets}
In this section we examine three modular invariants constructed by
A. N. Schelleken in \cite{Sch} which are based on level-rank
duality. These are entries 18, 27, and 40 in the table of
\cite{Sch}. Our goal in this section is to show that they can be
realized by conformal nets as an application of mirror extensions in
section \ref{mirrorextension}. For simplicity in this section we
will use $G_k$ to denote the corresponding conformal net $\A_{G_k}$
when no confusion arises.

\par
\subsection{Three mirror extensions}
\subsubsection{$\widetilde{SU(10)_2}$}\label{m1}
$\widetilde{SU(10)_2}$ is the simplest nontrivial example of mirror
extensions applying to $SU(2)_{10}\subset Spin(5)_1$ and
$SU(2)_{10}\times SU(10)_2\subset SU(20)_1$ in Theorem
\ref{mainmirror}. By Cor. \ref{mirrorindex} and Lemma \ref{lr}
$$
\mu_{\widetilde{SU(10)_2}}=20.
$$
Consider the induction for $SU(10)_2\subset \widetilde{SU(10)_2}.$
By Th. 5.7 of \cite{BEK} the matrix $Z_{\la\mu}=\lan \a_\la,
\a_\mu^{-}\ran$ commutes with the $S,T$ matrix of $SU(10)_2.$ Such
matrices are classified in \cite{Gan}, and it follows that there are
$15$ irreducible representations of $\widetilde{SU(10)_2}$ given as
follows: $\a_{J}^i, 0\leq i\leq 9, \a_{J}^i\sigma, 0\leq j\leq 4.$
The fusion rules are determined by the following relations:
$$
[\bar\sigma]= [\a_{J}^2 \sigma], [\a_{J}^5 \sigma]=[\sigma],
[\sigma\bar{\sigma}]= [1]+[\a_J^5]
$$
The restrictions of these representations to  $SU(10)_2$ are given
as follows:
$$
[\a_{J}^i\res \ ]= [J^i(2\La_0)]+[J^i(\La_3+\La_7)], 0\leq i\leq 9;
[\a_{J}^i\sigma]=[J^i(\La_0+\La_3)]+[J^i(\La_5+\La_8)] ,0\leq j\leq
4.
$$
It follows that modulo integers the conformal dimensions are given
as
$$
h_{\a_{J}^i}= \frac{i(10-i)}{10}, 0\leq i\leq 9 , h_\sigma=
\frac{77}{80},  h_{a_J\sigma}= \frac{25}{16},  h_{a_J^2\sigma}=
\frac{157}{80}, h_{a_J^3\sigma}= \frac{173}{80}=h_{a_J^4\sigma} .
$$
\begin{remark}
The modular tensor category (cf. \cite{Tu}) from representations of
$\widetilde{SU(10)_2}$ as given above seems to be unknown before. It
will be interesting to understand our construction from a
categorical point of view.
\end{remark}
The following simple lemma will be used later:
\begin{lemma}\label{so}
$\A_{Spin(n)_1}$ is a completely rational net whose irreducible
representations are in one to one correspondence with irreducible
representations of $\L Spin(n)_1.$ When $n$ is odd there are three
irreducible representations $1, \mu_0,\mu_1$ with index $1,1,
\sqrt{2}$ respectively and fusion rules $[\mu_1^2]=[1]+[\mu_0];$
when $n=4k+2, k\in {\mathbb N}$ the fusion rule is $\Z_4; $ when
$n=4k, k\in {\mathbb N}$ the fusion rule is  $\Z_2\times \Z_2.$
\end{lemma}
\proof By Th. 3.10 of \cite{Bo} it is enough to prove that
$\mu_{\A_{Spin(n)_1}} =4.$

When $n=5$ this follows from conformal inclusion $SU(2)_{10}\subset
Spin(5)_1$ and Lemma \ref{indexab}. Consider the inclusion
$SO(n)\times U(1)\subset SO(n+2).$ Note that the fundamental group
of $SO(n)$ is $\Z_2.$ It follows that loops with even winding
numbers in $LU(1)$ can be lifted to $LSpin(n),$ and we have a
conformal inclusion $LSpin(n-2)_1\times LU(1)_4\subset LSpin(n)_1.$
Since $\mu_{\A_{U(1)_4}}=4$ by \S 3 of \cite{X3m},  and the index of
$\A_{Spin(n-2)_1}\times \A_{U(1)_4}\subset \A_{Spin(n)_1}$ is
checked to be $2$, by induction one can easily prove the lemma for
all odd $n$. When $n$ is even we use the conformal inclusion
$\A_{SU(n/2)_1}\times \A_{U(1)_{2n}}\subset \A_{Spin(n)_1}$ with
index $n/2$. Note that $\mu_{\A_{SU(n/2)_1}}= n/2,
\mu_{\A_{U(1)_{2n}}}=2n$ by \S 3 of \cite{X3m}, and by Lemma
\ref{indexab} we have $\mu_{\A_{Spin(n)_1}} =4.$
\endproof

\subsubsection{$\widetilde{SU(9)_3}$}\label{m2}
$\widetilde{SU(9)_3}$ is an extension of $SU(9)_3$ by applying Th.
\ref{mainmirror} to $SU(3)_9\subset (E_6)_1$ and $SU(3)_9\times
SU(9)_3\subset SU(27)_1.$ By Cor. \ref{mirrorindex} and Lemma
\ref{lr} $ \mu_{\widetilde{SU(9)_3}}= 9. $ Recall the branching
rules for $SU(3)_9\subset (E_6)_1$ (We use $1_0$ to denote the
vacuum representation of $(E_6)_1$ and $1_+, 1_-$ the other two
irreducible representations of $(E_6)_1$):
$$
[1_0\res \ ]= \sum_{0\leq i\leq 2} ([\dot{J}^i(9\La_0)]
+[\dot{J}^i(\La_0+4\La_1+4\La_2)], \   [1_+\res \ ]=[1_-\res \ ]
\sum_{0\leq i\leq 2} ([\dot{J}^i(5\La_0+2\La_1+2\La_2)]
$$
where $\dot J:=9\La_1.$

Consider inductions with respect to
$$
\widetilde{SU(9)_3}\subset SU(9)_3 .$$ By Th. \ref{mainmirror} and
Lemma \ref{lr} the vacuum of $\widetilde{SU(9)_3}$ restricts to
representation
$$ \sum_{0\leq i\leq 2} ([{J}^{3i}(9\La_0)] +[J^{3i}(\La_3+\La_7+\La_8)]
$$ of $SU(9)_3.$  Since $J$ is local with  the above representation, by
Lemma \ref{a=a'} $\a_J$ is a DHR representation of
$\widetilde{SU(9)_3},$ and $[\a_J^3]=[1].$ One can determine the
remaining irreducible representations of $\widetilde{SU(9)_3}$ by
using \cite{Gan} as in \S\ref{m1}. Here we give a different approach
which will be useful in \S\ref{m3}. We note that $M(\dot{J})= J^3,
 M(\dot {J}^i(5\La_0+2\La_1+2\La_2))=
J^{3i}(\Lambda_4+\la_6+\Lambda_8), i=0,1,2$ by Lemma \ref{lr}, where
$M$ is defined as before Lemma \ref{M}.  By Prop. \ref{mirrora} we
have
$$
\lan \a_{\Lambda_4+\La_6+\Lambda_8},
\a_{\Lambda_4+\Lambda_6+\Lambda_8}^{-}\ran= 2.
$$
It follows that there are two irreducible DHR representations
$\tau_1,\tau_2$ of $\widetilde{SU(9)_3}$ such that
$\a_{\Lambda_4+\La_6+\Lambda_8}\succ [\tau_1]+[\tau_2],$ and
$\tau_1,\tau_2$ are the only two irreducible subsectors of
$\a_{\Lambda_4+\La_6+\Lambda_8}$ which are DHR representations. We
have for $i=1,2$ $\lan \tau_i,  \a_\mu^{-}\ran\leq \lan
\a_{\Lambda_4+\La_6+\Lambda_8},  \a_\mu^- \ran.$ Note that if the
color of $\mu$ is nonzero, then $\lan
\a_{\Lambda_4+\Lambda_6+\Lambda_8}, \a_\mu^-\ran=0$ by Lemma
\ref{3.3} since $\Lambda_4+\Lambda_6+\Lambda_8$ has color $0.$ If
$\mu$ has color $0$, by Lemma \ref{lr} and Prop. \ref{mirrora} we
have
$$
\lan \a_{\Lambda_4+\Lambda_6+\Lambda_8},  \a_\mu^-\ran
$$ is nonzero only when $\mu=J^{3i}(\Lambda_4+\Lambda_6+\Lambda_8), i=0,1,2. $
 It follows that
$$
\lan \tau_i, \a_\mu \ran =1
$$ when $\mu=J^{3i}(\Lambda_4+\Lambda_6+\Lambda_8), i=0,1,2,$ and
$$
\lan \tau_i, \a_\mu \ran =0
$$ when $\mu\neq J^{3i}(\Lambda_4+\Lambda_6+\Lambda_8), i=0,1,2,$
Hence the restriction of $\tau_i$ to $SU(9)_3$ are given as follows:
$$
[\tau_i\res \ ]=\sum_{0\leq j\leq
2}[J^{3j}(\Lambda_4+\Lambda_6+\Lambda_8)]
$$
It follows that the index of $\tau_i,i=1,2 $ is one, and since
$$
[(\a_J \tau_i)\res \ ]=\sum_{0\leq j\leq
2}[J^{3j+1}(\Lambda_4+\Lambda_6+\Lambda_8)],
$$
it follows that $[\a_j\tau_i]\neq [\tau_i].$ Hence the irreducible
representations of $\widetilde{SU(9)_3}$ are given by
$$ 1, \a_J, \a_J^2, \a_J^i \tau_k, 0\leq i\leq 2, k=1,2.$$
These representations generate an abelian group of order $9,$ it
must be either $\Z_3\times \Z_3$ or $\Z_9$. Note that by Lemma
\ref{3.3}
$$\lan \a_J, \tau_i^k\ran \leq \lan \a_J,
\a_{\Lambda_4+\Lambda_6+\Lambda_8}^k\ran =0, \forall k\geq 0
$$ since $J$ has color $3$ while $\Lambda_4+\Lambda_6+\Lambda_8$ has
color $0,$ it follows that these representations generate an abelian
group $\Z_3\times \Z_3.$  Modulo integers the conformal dimensions
of $\tau_k, \a_J$ are given by
$$
h_{\a_J}= \frac{4}{3}, h_{\tau_k}= \frac{7}{3},  h_{\a_J^2}=
\frac{7}{3}, h_{\a_J\tau_k}= \frac{11}{3}, h_{\a_J^2\tau_k}=
\frac{14}{3}, k=1,2.
$$
\subsubsection{$\widetilde{SU(8)_4}$}\label{m3}

From conformal inclusion $Spin(6)_8\subset Spin(20)_1$ and
$Spin(6)\simeq SU(4)$ we obtain conformal inclusion $SU(4)_8\subset
Spin(20)_1.$ For simplicity we use $(0), (5/4)_1, (5/4)_2, (1/2)$ to
denote irreducible representations of $Spin(20)_1$ with conformal
dimensions $0, 5/4,5/4, 1/2$ respectively. By comparing conformal
dimensions the branching rules for $SU(4)_8\subset Spin(20)_1$ are
given by:
\begin{align*}
[(0)\res \ ]& =\sum_{0\leq i\leq 3}
([\dot{J}^i]+[\dot{J}^i(4\La_0+\La_1+2\La_2+\La_3)]),\\
[(5/4)_1\res \ ]&=[(5/4)_2\res \ ]=\sum_{0\leq i\leq 3}
[\dot{J}^i(3\La_0+\La_1+2\La_2+3\La_3)], \\
[(1/2)\res \ ]&=\sum_{0\leq i\leq 3}
([\dot{J}^i(6\La_0+2\La_2)]+[\dot{J}^i(3\La_0+3\La_2+2\La_3)]).
\end{align*}
Note that all representations appearing above have color $0.$\par
$\widetilde{SU(8)_4}$ is the extension of $SU(8)_4$ by applying Th.
\ref{mainmirror} to  $SU(4)_8\subset Spin(20)_1$ and $SU(4)_8\times
SU(8)_4\subset SU(32)_1.$ By Lemma \ref{lr} the spectrum of
$SU(8)_4\subset \widetilde{SU(8)_4}$ is given by
$$
\sum_{0\leq i\leq 3}
([{J}^{2i}]+[{J}^{2i}(\La_0+\La_4+2\La_5+\La_7)]).
$$
By Lemma \ref{so} and Lemma \ref{indexab}
$\mu_{\widetilde{SU(8)_4}}= 8.$

By using Prop. \ref{mirrora} similar as in \S\ref{m2} we obtain all
irreducible representations of $\widetilde{SU(8)_4}$ as follows:
$$
1, \a_J, (3/4)_1, (3/4)_2, (1/2), \a_J (3/4)_1, \a_J (3/4)_2, \a_J
(1/2).
$$
%which generate an abelian group of order $8$ under compositions.
These representations restrict to $SU(8)_4$ as follows:
\begin{align*}
[\a_J\res \ ]& =\sum_{0\leq i\leq 3}
([{J}^{2i+1}]+[{J}^{2i}(\La_0+\La_4+2\La_5+\La_7)]),\\
[(\a_J^j(3/4)_k)\res \ ]& =\sum_{0\leq i\leq 3}
[{J}^{2i}(\La_0+\La_3+2\La_6+\La_7)]), j=0,1, k=1,2;\\
[(\a_J^j(1/2))\res \ ]& =\sum_{0\leq i\leq 3} ([{J}^{2i+j}(2\La_0+
\La_3+\La_5)]+[{J}^{2i+j}(2\La_5+2\La_7)]), j=0,1.
\end{align*}
The conformal dimensions modulo integers are as follows:
$$
h_{\a_J}= 7/4, h_{(3/4)_k}= 3/4, k=1,2, h_{(1/2)}=1/2,
h_{\a_J(3/4)_1}= h_{\a_J(3/4)_2}= 5/2, h_{\a_J(1/2)}=9/4,
$$ which explain our notations.
The irreducible representations of $\widetilde{SU(8)_4}$ generate an
abelian group of order $8$ under compositions, so the abelian group
is $\Z_2\times \Z_2\times \Z_2, \Z_2\times \Z_4 $ or $\Z_8.$  By
Lemma \ref{3.3} $\lan \a_J, (3/4)_k^j\ran = \lan \a_J, (1/2)^j\ran=
0, k=1,2,\forall j\geq 0$ since the restriction of $\a_J$ to
$SU(8)_4$ has color $4$ while the restriction of $(3/4)_k, (1/2)$ to
$SU(8)_4$ has color $0$, it follows that  $\Z_8$ is impossible. Note
that the conjugate of $(1/2)$ has conformal dimension $1/2$, and it
must be $(1/2)$, so $[(1/2)^2]= [1].$ To rule out the possibility of
$\Z_2\times \Z_4 ,$ note that this can only happen when the order of
$(3/4)_1$ is $4$, and we must have $[(1/2)]=[(3/4)_1^2],
[(3/4)_2]=[(3/4)_1^3].$  By monodromy equation we have
$$
\e((3/4)_1,(3/4)_1)^2= 1, \e((3/4)_1, (1/2))\e((1/2),(3/4)_1)=-1.
$$
On the other hand by Lemma 4.4 of \cite{Rehren} we have
$$
\e((3/4)_1, (1/2))\e((1/2),(3/4)_1)= \e((3/4)_1,
(3/4)_1^2))\e((3/4)_1^2,(3/4)_1)=\e((3/4)_1,(3/4)_1)^4=1,
$$
a contradiction. It follows that irreducible representations of
$\widetilde{SU(8)_4}$ generate $\Z_2\times \Z_2\times \Z_2$ under
compositions, and we have
$$
\overline{[(3/4)_1]}= [(3/4)_1], [(1/2)(3/4)_1]=[(3/4)_2].
$$
\subsection{Further extensions by simple currents}
\subsubsection{No. 40 of \cite{Sch}}\label{40}
The modular invariant No. 40 in \cite{Sch} suggests that we look for
simple current extensions of $\widetilde{SU(9)_2}\times
SU(5)_1\times SO(7)_1.$ For simplicity we use $y^i=\La_i, 0\leq
i\leq 9 $ to denote the irreducible representation of $SU(5)_1.$
Note that $h_{y^2} =3/5.$ We use $(1/2), (7/16)$ to denote the
irreducible representations of $SO(7)_1$ with conformal dimensions
$1/2, 7/16.$ Note that the index of $(1/2), (7/16)$ are $1, 2$
respectively. By \S\ref{m1} the conformal dimension of $u=(\a_J,
y^2, (1/2))$ is $h_{\a_J}+ h_y+ 1/2= 2.$ It follows that $u^i, 0\leq
i\leq 9$ is a local system of automorphisms. By Prop. \ref{simple}
there is a M\"{o}bius extension $\D=(\widetilde{SU(9)_2}\times
SU(5)_1\times SO(7)_1)\ltimes \Z_{10}$ of $\widetilde{SU(9)_2}\times
SU(5)_1\times SO(7)_1.$ By Cor. \ref{mirrorindex} and Lemma \ref{so}
$\mu_\D= 4.$ Consider now the inductions for
$$
\widetilde{SU(9)_2}\times SU(5)_1\times SO(7)_1\subset \D
$$
By using formulas for conformal dimensions in \S\ref{m1} one checks
easily that
$$H((\sigma,y^3,(7/16)), u) =H((1,1,(1/2)),u)=1.$$
By Lemma \ref{a=a'} we conclude that $\a_{(\sigma,y^3,(7/16))},
\a_{(1, 1, (1/2))}$ are DHR representations of $\D$ with index $2,
1$ respectively. Note that by Lemma \ref{3.3}
$$
\lan \a_{(\sigma,y^3,(7/16))}, \a_{(\sigma,y^3,(7/16))}\ran
=\sum_{0\leq i\leq 9}\lan (\sigma,y^3,(7/16)),
(\sigma,y^3,(7/16))u^i\ran = 2
$$
where in the last step we have used $[\sigma a_J^5]=[1].$ It follows
that $[\a_{(\sigma,y^3,(7/16))}]= [\delta_1]+[\delta_2].$ Since
$\mu_\D=4,$ the list of irreducible representations are given by
$$
1,\a_{(1, 1, (1/2))}, \delta_1, \delta_2.
$$
The conformal dimensions modulo integers are $h_{\delta_1}=
h_{\delta_2}=1, h_{\a_{(1, 1, (1/2))}}=1/2.$ These representations
generate an abelian group of order $4$. To rule out $\Z_4$, note
that $[\a_{(1, 1, (1/2))}^2]=[1].$ Without losing generality we
assume that $\delta_1$ has order $4$. Then we must have
$[\delta_1^2]=[\a_{(1, 1, (1/2))}], [\delta_1^3]=[\delta_2]. $ By
monodromy equation we have $\e(\delta_1,\delta_1)^2= -1,
\e(\delta_1, \delta_2)\e(\delta_2, \delta_1)=1.$ On the other hand
by Lemma 4.4 of \cite{Rehren} we have $ \e(\delta_1,
\delta_2)\e(\delta_2, \delta_1)=\e(\delta_1,\delta_1)^6 = -1,$ a
contradiction. In particular we have $[\delta_1^2]=[1].$ Hence $1,
\delta_1$ is a local system of automorphisms, and by Prop.
\ref{simple} we conclude that the there is further extension
$\D\ltimes \Z_2$ of $\D$. By Lemma \ref{indexab} we have
$\mu_{\D\ltimes \Z_2}=1,$ i.e., $\D\ltimes \Z_2$ is holomorphic. The
spectrum of $SU(10)_2\times SU(5)_1\times Spin(7)_1\subset \D\ltimes
\Z_2$ is given by entry 40 in the table of \cite{Sch}:
\begin{equation*}
\sum_{o\leq i\leq 9} ([(J^i, y^{2i}, (1/2)^i)]+[(J^i(\La_3+\La_7),
y^{2i}, (1/2)^i)]+ [(J^i(\La_3+\La_6), y^{2i+4}, (7/16))])
\end{equation*}
\subsubsection{No. 27 of \cite{Sch}}\label{27}
No. 27 in the table of \cite{Sch} suggests that we look for simple
current extensions of $\widetilde{SU(9)_3}\times SU(3)_1\times
SU(3)_1.$ Label irreducible representations of $SU(3)_1$ by their
conformal dimensions as $1, (1/3)_1, (1/3)_2.$ Denote by $x_1=(\a_J,
(1/3)_1, (1/3)_1), x_2=(\tau_1, (1/3)_1, (1/3)_2.$  By using
formulas for conformal dimensions in \S\ref{m2} and Lemma
\ref{checklocal} it is easy to check that the following set $ x_1^i
x_2^j, 0\leq i,j\leq 2$ is a local system of automorphisms. Hence by
Prop. \ref{simple}  there is a M\"{o}bius extension
$\D_1=(\widetilde{SU(9)_3}\times SU(3)_1\times SU(3)_1) \ltimes
(\Z_3\times\Z_3)$ of $\widetilde{SU(9)_3}\times SU(3)_1\times
SU(3)_1$ with spectrum $\sum_{0\leq i,j\leq 2}[x_1^i x_2^j].$ By
Lemma \ref{indexab} $\mu_{\D_1}= 1,$ so $\D_1$ is holomorphic. The
spectrum of ${SU(9)_3}\times SU(3)_1\times SU(3)_1\subset \D_1$ is
given by (entry (27) of \cite{Sch}):
\begin{align*}
\sum_{ 0\leq i\leq 9} & ([(J^i, (1/3)_1^i, (1/3)_1^i)] + [(J^i(\La_4
+\La_6+\La_8), (1/3)_1^{i-1}, (1/3)_1^{i+1})] \\
&+ [(J^i(\La_4 +\La_6+\La_8), (1/3)_1^{i+1}, (1/3)_1^{i-1})] +
[(J^i(\La_3+\La_7+\La_8), (1/3)_1^i, (1/3)_1^i)]
\end{align*}
\begin{remark}
One can choose other local systems of automorphisms which generate
$\Z_3\times \Z_3.$  For an example one such choice is a  local
system of automorphisms  given by $x_1'^i {x_2'}^j,0\leq i,j\leq 2$
with $x_1'=(\a_J, (1/3)_1, (1/3)_2), x_2'=(\tau_1, (1/3)_1,
(1/3)_1).$ However by remark \ref{outer} it is easy to check that
the corresponding extension is simply $\Ad U(\D_1)$ which is
isomorphic to $\D_1$, where $\Ad U$ implements the outer
automorphism of the last factor of $SU(3)_1.$ Similar statement
holds for other choices of local systems of automorphisms which
generate $\Z_3\times \Z_3.$
\end{remark}
\subsubsection{No. 18 of \cite{Sch}}\label{18}
No. 18 in the table of \cite{Sch} suggests that we look for simple
current extensions of $\widetilde{SU(8)_4}\times SU(2)_1\times
SU(2)_1\times SU(2)_1.$ As before we label the non-vacuum
representation $(1/4)$ of $SU(2)_1$ by its conformal dimension. Set
$z_1=(\a_J, (1/4),0,0), z_2=((3/4)_1, 0, (1/4),0), z_3= ((3/4)_2,
0,0, (1/4)).$ Then by the formulas for conformal dimensions and
fusion rules in \S\ref{m3} one checks easily that $H( z_i, z_j)=1,
1\leq i,j\leq 3.$ Hence $\{ z_1,z_2,z_3\}$ generate an abelian group
$\Z_2\times \Z_2\times \Z_2$ which is a local system of
automorphisms by Lemma \ref{checklocal}. By Prop. \ref{simple} we
conclude that there is a M\"{o}bius extension
$\D_2:=(\widetilde{SU(8)_4}\times SU(2)_1\times SU(2)_1\times
SU(2)_1)\ltimes (\Z_2\times \Z_2\times \Z_2).$ By Lemma
\ref{indexab} we have $\mu_{\D_2}=1,$ i.e., $\D_2$ is holomorphic.
The spectrum of ${SU(8)_4}\times SU(2)_1\times SU(2)_1\times
SU(2)_1\subset \D_2$ is given by (entry (18) of \cite{Sch}):
\begin{align*}
\sum_{0\leq i\leq 7}  ([(J^i, (1/4)^i,0,0)]+
[(J^i(\La_0+\La_4+\La_5+\La_7), (1/4)^i,0,0)]
\\
+[(J^i(\La_5+2\La_7), J_1^i,(1/4),(1/4))]
+[(J^i(2\La_0+\La_3+\La_5), (1/4)^i,(1/4),(1/4))] \\
 + [(J^i(\La_0+\La_3+\La_6+\La_7), (1/4)^i,0,(1/4))] +
[(J^i(\La_0+\La_3+\La_6+\La_7), (1/4)^i, (1/4),0)])
\end{align*}
\subsubsection{The main Theorem}
By Lemma \ref{extconformal} $\D\ltimes \Z_2, \D_1, \D_2$ as
constructed in \S\ref{40}, \S\ref{27} and \S\ref{18} are in fact
conformal nets since they contain conformal subnets with finite
index, and in summary we have proved the following:
\begin{theorem}\label{main}
There are holomorphic conformal nets (with central charge $24$)
which are conformal extensions of $SU(10)_2\times SU(5)_1\times
Spin(7)_1, SU(9)_3\times SU(2)_1\times SU(2)_1, {SU(8)_4}\times
SU(2)_1\times SU(2)_1\times SU(2)_1$ with spectrum given by the
representations at the end of \S\ref{40}, \S\ref{27} and \S\ref{18}
respectively.
\end{theorem}
\subsection{Two conjectures}
The holomorphic conformal net corresponding to $V^\sharp$ of
\cite{FLM} was constructed in \cite{KL2}. This net can also be
constructed using the result of \cite{DX} as a simple current $\Z_2$
extension of a $\Z_2$ orbifold conformal net associated with Leech
lattice given in \cite{DX}. Our first conjecture is an analogue of
the conjecture in \cite{FLM} for $V^\sharp$:
\begin{conjecture}\label{conj1}
Up to isomorphism there exists a unique holomorphic conformal net
with central charge $24$ and no elements of weight one.
\end{conjecture}
Our second conjecture is motivated by the results of \cite{Sch}:
\begin{conjecture}\label{conj2}
Up to isomorphism there exists  finitely many   holomorphic
conformal nets with central charge $24$ .
\end{conjecture}
Note that if one can obtain a theorem like the theorem in \S2 of
\cite{Sch} in the setting of conformal nets, then modulo conjecture
(\ref{conj1}) conjecture ({\ref{conj2}) is reduced to show that up
to equivalence, there are only finitely many conformal extensions of
a given completely rational net, and this should be true in view of
the results of \cite{IK}. However new methods have to be developed
to carry though this idea.

\end{document}